\documentclass[12pt,a4,amsmath, amsthm, amssymb]{amsart}

\newtheorem{Lemma}{Lemma}[section]
\newtheorem{Theorem}[Lemma]{Theorem}
\newtheorem{Proposition}[Lemma]{Proposition}
\newtheorem{Corollary}[Lemma]{Corollary}

\newtheorem{remark}[Lemma]{Remark}
\newtheorem{definition}[Lemma]{Definition}
\newtheorem{example}[Lemma]{Example}

\newtheorem{Fact}[Lemma]{Fact}

\def\bt{\begin{Theorem}}
\def\et{\end{Theorem}}
\def\bl{\begin{Lemma}}
\def\el{\end{Lemma}}
\def\bp{\begin{Proposition}}
\def\ep{\end{Proposition}}
\def\bcor{\begin{Corollary}}
\def\ecor{\end{Corollary}}
\def\bpf{\begin{proof}}
\def\epf{\end{proof}}

\def\brem{\begin{remark}}
\def\erem{\end{remark}}

\def\bedef{\begin{definition}}
\def\endef{\end{definition}}

\def\beg{\begin{example}}
\def\eeg{\end{example}}

\def\bef{\begin{Fact}}
\def\eef{\end{Fact}}

\def\bc{\begin{center}}
\def\ec{\end{center}}
\def\noi{\noindent}

\def\vsq{\vskip .25cm}
\def\beq{\begin{equation}}
\def\eeq{\end{equation}}
\def\beqarray{\begin{eqnarray*}}
\def\eeqarray{\end{eqnarray*}}
\def\<{\leftangle}
\def\>{\rightangle}
\def\({\left(}
\def\){\right)}
\def\f{\varphi}

\def\<{\langle}
\def\>{\rangle}
\def\q{\quad}
\def\r{\rho}
\def\a{\alpha}
\def\b{\beta}

\def\p{\pi}
\def\d{\delta}
\def\h{\hbox}

\def\t{\tau}
\def\r{\rho}
\def\l{\lambda}

\def\w.r.t.{with respect to}
\def\R{{\mathbb{R}}}
\def\N{{\mathbb{N}}}

\def\C{{\mathbb{C}}}

\def\bq{\begin{quote}}
\def\eq{\end{quote}}

\def\bit{\begin{itemize}}
\def\eit{\end{itemize}}
\def\i{\item}
\def\ben{\begin{enumerate}}
\def\een{\end{enumerate}}

\usepackage{float}
\usepackage{epsfig}
\begin{document}

\title[Time-fractional backward heat conduction problem]
{A new regularisation for time-fractional backward heat conduction problem}
\author{M.Thamban Nair and P. Danumjaya}
\address{Department of Mathematics, BITS Pilani, K.K. Birla Goa Campus, Zuarinagar, Goa 403726, INDIA}
\email{mtnair@goa.bits-pilani.ac.in; danu@goa.bits-pilani.ac.in}
\today

\begin{abstract} 
It is well-known that the backward heat conduction problem of recovering the temperature $u(\cdot, t)$ at a time $t\geq 0$ from the knowledge of the temperature at a later time, namely $g:= u(\cdot, \t)$ for  $\t>t$,   is ill-posed, in the sense that small error in $g$ can lead to large deviation in $u(\cdot, t)$. However, in the case of a time fractional backward heat conduction problem (TFBHCP), the above problem is well-posed for $t>0$ and ill-posed for $t=0$. We use this observation to obtain  stable approximate solutions for the TFBHCP for $t=0$, and derive error estimates under suitable source conditions. We shall also provide some numerical examples to illustrate the approximation properties of the  regularized solutions.   
\end{abstract}
\maketitle



\section{Introduction}

For $0<\a<1$, consider the time-fractional heat equation 
\beq\label{eq-1} \frac{\partial ^\a u}{\partial t^\a} =  \frac{\partial^2 u}{\partial x^2},\q 0<x<\pi.
\eeq
In the above, we used the $\a$-derivative of $u$ \w.r.t. $t$ in the  {\it Caputo sense}. That is, if $\f$ is a real valued differentiable  function on an open interval of the form $(0, a)$ for some $a>0$,   
$$\frac{d^\a\f}{dt^\a} (t) = \frac{1}{\Gamma(1-\a)}\int_0^t \frac{\f'(s)}{(t-s)^\a} ds,\q 0<t<a.$$ 
It is to be observed that for $\a=1$, the equation (\ref{eq-1}) reduces to the ordinary heat equation
$$ \frac{\partial  u}{\partial t} =  \frac{\partial^2 u}{\partial x^2},\q 0<x<\pi,$$ 
and in that case, under the boundary condition 
\beq\label{BC} u(0,t) = 0 = u(\pi, t),\q t>0,\eeq
and initial condition  
\beq\label{IC} u(0, x) = f_0(x),\q 0<x<\pi\eeq
for $f_0\in L^2[0, \pi]$, 
the  solution $u(\cdot, \cdot)$ is given by 
\beq\label{Heat}u(\cdot,t) = \sum_{n=1}^\infty E(-\l_n^2t) \<f_0, \f_n\>\f_n,\eeq
where $E(\cdot)$ is the exponential function $E(s):=e^s$ and 
$$\l_n:= n,\q \f_n(x) = \sqrt{\frac{2}{\pi}}\sin (kx).$$
We know that $\{\f_n: n\in \N\}$ is an orthonormal basis of $L^2[0, \pi]$. 

It is known (cf.\,\cite{Hao})  that, for the time-fractional heat equation (\ref{eq-1}) along with the boundary and initial conditions (\ref{BC})-(\ref{IC}), similar representation for its solution is possible by replacing $E(-\l_n^2t)$ in  (\ref{Heat}) by $E_\a(-\l_n^2t^\a)$, where  $E_\a(\cdot)$ is   the {\it Mittag-Leffler function} defined by
\beq\label{ML} E_\a(z) = \sum_{k=0}^\infty  \frac{z^k}{\Gamma(k\a +1)},\q z\in \C,\eeq 
with $\Gamma(\cdot)$ being the gamma function 
$$\Gamma(s):=\int_0^\infty e^{-t}t^{s-1}dt,\q s>0.$$
Clearly, $E_1(z) = E(z)$ for $z\in \C$. In fact,  for any $f_0\in L^2[0, \pi]$,  the solution  $u_\a(\cdot, \cdot)$ of (\ref{eq-1}) along with the boundary and initial conditions (\ref{BC})-(\ref{IC}), is given by 
\beq\label{eq-2} u_\a(\cdot,t) = \sum_{n=1}^\infty E_\a(-\l_n^2t^\a) \<f_0, \f_n\>\f_n\eeq
with  $\<\cdot, \cdot\>$ being the $L^2$-inner product (cf.\,\cite{Hao, Tuan-1}). Throughout, by $\|\cdot\|$ we shall mean the norm induced by the $L^2$-inner product.

With regard to the Mittag-Leffler function defined in (\ref{ML}), we have the following result which we shall make use of throughout the paper. 

\bl\label{lem-1} {\rm(cf. \cite{Liu-Yama, FDE-Pod})}
Given real numbers $\a_0,\a_1$ such that  $0<\a_0<\a_1<1$, there exist $C_1>0$ and $C_2>0$ such that 
$$ \frac{C_1}{\Gamma(1-\a)(1+x)} \leq E_\a(-x) \leq   \frac{C_2}{\Gamma(1-\a)(1+x)}$$
for all $x>0$ and for all $\a\in [\a_0, \a_1]$. 
\el 
In view of the above lemma, it can be seen that $u_\a(\cdot, t)$  in (\ref{eq-2} ) is a $C^\infty$-function for each $t>0$. 

Analogous to the well-studied  ordinary {\it backward heat conduction problem}, let us consider the  following inverse  problem,   the  {\it time-fractional backward heat conduction problem}  (TFBHCP) associated with (\ref{eq-1}) and (\ref{BC}):

\vsq
\noi
{\bf Problem ($P_t$):} Knowing $g:=u_\a(\cdot, \t)$ for some $\t>0$, find $ u_\a(\cdot, t)$ for $0\leq t<\t$. 
\vsq
Many studies have shown that fractional diffusion equation model is appropriate for investigating problems arising in the areas of  spatially disordered systems, porous media, fractal media, turbulent fluids and plasmas, biological media with traps, and  stock price movements and so on (see \cite{FDE-Pod, Tuan-1}, and the references therein).  The regularization theory for the inverse problems associated with fractional-order PDEs is still in its infancy. 

We shall see that the inverse problem $P_t$  is well-posed if $0<t<\t$ and ill-posed if $t=0$. This observation and the subsequent analysis lead us to the conclusion that $\{P_t: 0<t<\t\}$ gives a  {\it regularization family} for  obtaining  stable approximate solutions for the ill-posed inverse problem $(P_0)$. To our knowledge, no study is carried out using the above observation, though various regularization methods are discussed recently (see, e.g. \cite{Tuan-1, {Tuan-Dang}, {Hao}, kokila-nair}, and the references there in). We shall also provide estimates for the error 
$\|u_\a(t,\cdot) - u_\a(0, \cdot)\|_2$ under certain a priori source condition.


An outline of this paper is as follows. In section 2, we discuss the ill-posedness of the time fractional backward heat conduction problem (TFBHCP). Section 3 deals an operator theoretic formulation of the inverse problems. The new regularization family for the ill-posed inverse problem $P_0$ is introduced and its convergence is proved in  Section 4. In Section 5 and 6, we derive the error estimates for the noisy data and source conditions, respectively. Finally, we perform some numerical experiments to validate the theoretical results in Section 7.

\section{Ill-Posednss of the Inverse Problem}

Let $g=u_\a(\cdot, \t)$. Then, from equation (\ref{eq-2}), we have 
\beq\label{eq-3} g =  \sum_{n=1}^\infty E_\a(-\l_n^2\t^\a) \<f_0, \f_n\>\f_n.\eeq
Hence, 
\beq\label{eq-4} \<f_0, \f_n\> = \frac{\< g, \f_n\> } {E_\a(-\l_n^2\t^\a)} \q\forall \, n\in \N.\eeq
This shows that  $g := u_\a(\cdot, \t)$ must satisfy  the {\it Picard  condition}:
\beq\label{Picard}\sum_{n=1}^\infty \frac{|\< g, \f_n\>|^2 } {E_\a(-\l_n^2\t^\a)^2}<\infty.\eeq
By Lemma \ref{lem-1}, it is to be observed that 
$$ \frac{1}{C_2}   \Gamma (1-\a)(1+\l_n^2 \t^\a)  \leq \frac{1 } {E_\a(-\l_n^2\t^\a)^2} \leq \frac{1}{C_1}   \Gamma (1-\a)(1+\l_n^2 \t^\a).$$
Therefore, Picard  condition (\ref{Picard}) on $g$  is equivalent to the requirement 
\beq\label{Picard -1} \sum_{n=1}^\infty (1+\l_n^2\t^\a) |\< g, \f_n\>|^2 <\infty\eeq
which is again equivalent to 
$$\sum_{n=1}^\infty n^2  |\< g, \f_n\>|^2 <\infty.$$
Using (\ref{eq-4}), the representation of $u_\a(\cdot, t)$ in  (\ref{eq-2}) takes the form 
\beq\label{eq-5-0}  u_\a(\cdot,t) = \sum_{n=1}^\infty \frac{ E_\a(-\l_n^2t^\a) } {E_\a(-\l_n^2\t^\a)}  \< g, \f_n\> \f_n\eeq
so that 
\beq\label{eq-5} \|u_\a(\cdot,t)\|^2  = \sum_{n=1}^\infty  \Big|\frac{ E_\a(-\l_n^2t^\a) } {E_\a(-\l_n^2\t^\a)} \Big|^2| \< g, \f_n\>|^2.
\eeq
Again, using Lemma \ref{lem-1}, we have
\beq\label{eq-6} 
  \frac{C_1}{C_2}   \frac{(1+\l_n^2 \t^\a)}{ (1+\l_n^2 t^\a)} \leq  \frac{ E_\a(-\l_n^2t^\a) } {E_\a(-\l_n^2\t^\a)}  \leq  \frac{C_2}{C_1}   \frac{(1+\l_n^2 \t^\a)}{ (1+\l_n^2 t^\a)} 
 \eeq
Note that,  
$$ \frac{(1+\l_n^2 \t^\a)}{ (1+\l_n^2 t^\a)} \geq  \frac{\l_n^2 \t^\a}{ (\l_n^2+\l_n^2 t^\a)} = \frac{ \t^\a}{ (1+ t^\a)},$$
and for $0<t\leq \t$,
$$   \frac{(1+\l_n^2 \t^\a)}{ (1+\l_n^2 t^\a)}  \leq  \frac{(\l_n^2+\l_n^2 \t^\a)}{ \l_n^2 t^\a} =   \frac{(1+ \t^\a)}{  t^\a}.$$
Hence, for $0<t\leq \t$, 
\beq\label{eq-6-1} 
 \frac{C_1}{C_2}   \frac{ \t^\a}{ (1+ t^\a)}   \leq  \frac{ E_\a(-\l_n^2t^\a) } {E_\a(-\l_n^2\t^\a)}    \leq \frac{C_2}{C_1}   \frac{(1+ \t^\a)}{  t^\a}.
 \eeq
Hence, the representation (\ref{eq-5}) together with (\ref{eq-6-1}) imply  that if $0<t\leq \t$, then 
\beq\label{eq-7}  \frac{C_1}{C_2}   \frac{ \t^\a}{ (1+ t^\a)}   \|g\|  \leq \|u_\a(\cdot, t)\|   \leq  \frac{C_2}{C_1}   \frac{(1+ \t^\a)}{  t^\a} \|g\| \eeq
and if  $t=0$, then (\ref{eq-5}) and (\ref{eq-6}) imply 
\beq\label{eq-8} \|u_\a(\cdot,0)\|^2  \geq  \Big(\frac{C_1}{C_2} \Big)^2 \sum_{n=1}^\infty  (1+\l_n^2 \t^\a)^2 | \< g, \f_n\>|^2.\eeq
Note that, by (\ref{Picard -1}) the series on the right hand side of the above inequality converges. 
However, corresponding to an initial temperature  $\tilde u_\a(\cdot, 0)$, if the temperature at time $\t$ is  $\tilde g$, then the above arguments  lead to 
$$ \|u_\a(\cdot,0) - \tilde u_\a(\cdot,0)\|^2  \geq  \Big(\frac{C_1}{C_2} \Big)^2 \sum_{n=1}^\infty  (1+\l_n^2 \t^\a)^2 | \< g-\tilde g, \f_n\>|^2.$$
From this we obtain 
\beq\label{eq-9} \|u_\a(\cdot,0)- \tilde u_\a(\cdot,0) \|   \geq   (1+\l_n^2 \t^\a) \Big(\frac{C_1}{C_2}  \Big)\|g-\tilde g\|\q\forall\, n\in \N.\eeq
This shows that small error in $g$ can lead to large deviation in the solution $u_\a(\cdot,0)$, even when the the data satisfy the Picard condition  as in (\ref{Picard}). 
Thus, from the  inequalities in (\ref{eq-7}) and (\ref{eq-9}), we can infer the following.

\bt 
Let $C_1, C_2$ and $\a$ be as in Lemma \ref{lem-1}. Then  the TFBHCP   $(P_t)$ is well-posed for $0<t<\t$ and the problem $(P_0)$ is ill-posed.
\et 

What we are interested in  is to find stable approximate solutions for the ill-posed inverse problem $(P_0)$. 

\section{Operator Theoretic formulation of the inverse problems}  

For a few  observations on operators  that appear in this section,  we shall make use of the following proposition  based on basic results from functional analysis.  For the sake of completeness of the exposition, we provide its proof as well. 

\bp\label{prop-1}
Let ${\mathcal H}$ be an infinite  dimensional separable Hilbert space  and let $\{v_n: n\in \N\}$ be an orthonormal basis of ${\mathcal H}$. Let $(\mu_n)$ be a bounded sequence of  real numbers and $A: {\mathcal H}\to {\mathcal H}$ be defined by 
$$Av = \sum_{n=1}^\infty \mu_n \<v, v_n\>_{\mathcal H}v_n,\q v\in \mathcal H,$$
where  $\<\cdot, \cdot \>_{\mathcal H}$ denotes the inner product on ${\mathcal H}$. Then $A$ is a self-adjoint, bounded linear operator.  Further, we have the following:
\ben 
\i[\rm (i)] If $\mu_n\to 0$, then $A$ is a compact operator.
\i[\rm (ii)] If $\mu_n\not=0$ for all $n\in \N$, then $A$ is injective and range of $A$ is dense.
\i[\rm (iii)] If there exists $c_0>0$ such that $|\mu_n|\geq c_0$ for all $n\in \N$, then $A$ is injective,  range of $A$ is closed, and its inverse from the range is continuous.    
\een
In particular, if assumptions in (ii) and (iii) are satisfied, then $A$ is bijective and its inverse is continuous.
\ep 

\bpf 
Using the boundedness of $(\mu_n)$ it follows from Riesz-Fischer theorem (cf. \cite{Nair-fa}) that $A$ bounded linear operator on   ${\mathcal H}$ with $\|A\| \leq \sup\{|\mu_n|: n\in \N\}$. Also, since $\mu_n\in \R$, we that $A$ is a self-adjoint operator. 

(i) Suppose $\mu_n\to 0$. For $n\in \N$, let  $A_n: {\mathcal H}\to {\mathcal H}$ be defined by
$$A_nv = \sum_{j=1}^n \mu_j \<v, v_j\>_{\mathcal H}v_j,\q v\in \mathcal H.$$
Then we see that 
$$\|A-A_n\|\leq \sup\{|\mu_j|: j>n\}\to 0.$$ 
Since each $A_n$ is a finite rank operator, it follows that (cf. Theorem  9.-  in \cite{Nair-fa}) $A$ is a compact operator . 

(ii)  Suppose   $\mu_n\not =0$ for all $n\in \N$. Then for $v\in \mathcal H$,  we have 
$$Av=0 \iff \mu_n \<v, v_n\>_{\mathcal H} =0\,\,\forall \, n\in \N\iff \<v, v_n\>_{\mathcal H} =0\,\, \forall \, n\in \N\iff v=0.$$
Hence,  $A$ is injective.  Now, to see that $R(A)$, the range of $A$,  is dense in ${\mathcal H}$, let $w\in {\mathcal H}$ be  such that $\<Av, w\>=0$  for all $v\in {\mathcal H}$. Then, in particular, we have 
$$\mu_n\<\f_n, w\> = \<A\f_n, w\>=0\q\forall n\in \N.$$
From this, using again the fact that $\mu_n\not =0$ for all $n\in \N$, we have $w=0$. Thus, we have proved that $R(A)^\perp = \{0\}$, which implies, by projection theorem, that $R(A)$ is dense. 

(iii) Suppose there exists $c_0>0$ such that $|\mu_n|\geq c_0$ for all $n\in \N$.  Then we have
$$\|Av\|_{\mathcal H} \geq  c_0\|v\|_{\mathcal H}\q\forall\, v\in \mathcal H.$$
From this, the conclusions follow. 

The last part of the theorem is obvious. 
\epf 


Now, consider the operator $A_\a: L^2[0, \pi]\to L^2[0, \pi]$ defined by 
$$A_\a f = \sum_{n=1}^\infty E_\a(-\l_n^2\t^\a) \<f, \f_n\>\f_n,\ f\in L^2[0, \pi].$$
Using  Lemma \ref{lem-1}, we see that  
$$E_\a(-\l_n^2\t^\a)\to 0\q\h{as}\q n\to \infty.$$
Hence, by Proposition \ref{prop-1}, $A_\a$   
is a compact operator of infinite rank.  Hence, in view of equation (\ref{eq-3}), the fact that the problem $(P_0)$ is ill-posed also follows from the observation that (\ref{eq-3}) is same as solving the compact operator equation 
\beq\label{op-eq-0} A_\a f=g,\eeq
which is an ill-posed problem. 

Again, by Proposition \ref{prop-1},  $A_\a$ is one-one and its range is dense in $L^2[0, \pi]$. Hence, if  $g\in L^2[0, \pi]$ satisfies the the Picard condition (\ref{Picard}), then it  is in the range of $A_\a$ and  $f_0:=u_\a(\cdot, 0)$ is the the generalized solution of the operator equation (\ref{op-eq-0}), that is, 
$$f_0 = A_\a^\dagger g,$$
where $A_\a^\dagger$ denotes the Moore-Penrose inverse of $A_\a$ (cf. \cite{Nair-linop}). 

Next, we observe from (\ref{eq-2}) that  
$$\<f_0, \f_n\> = \frac{\<u_\a(\cdot, t), \f_n\>}{E_\a(\l_n^2t^\a)}\q\forall\, n\in \N,$$
so that  equation  (\ref{eq-3}) leads to  
\beq\label{eq-well-op} g= \sum_{n=1}^\infty \frac{E_\a(-\l_n^2\t^\a)}{E_\a(-\l_n^2t^\a)} \<u_\a(\cdot, t), \f_n\> \f_n. \eeq  
Interchanging  $t$ and $\t$ in (\ref{eq-6-1}), we obtain 
$$ 
\frac{C_1}{C_2}   \frac{ t^\a}{ (1+ \t^\a)}   \leq  \frac{ E_\a(-\l_n^2\t^\a) } {E_\a(-\l_n^2t^\a)}    \leq \frac{C_2}{C_1}   \frac{(1+ t^\a)}{  \t^\a}.$$
Hence, by Proposition \ref{prop-1}, it follows that  for $0<t< \t$,  $B_{\a,t}: L^2[0, \pi]\to L^2[0, \pi]$ defined by 
$$B_{\a,t} f = \sum_{n=1}^\infty \frac{E_\a(-\l_n^2\t^\a)}{E_\a(-\l_n^2t^\a)} \<f, \f_n\> \f_n, \q f\in L^2[0, \pi], $$
is a bijective bounded linear operator with continuous inverse.   Thus, the problem $(P_t)$ of recovering $u_\a(\cdot, t)$ from $g$, which corresponds to the equation (\ref{eq-well-op}) is same as the problem of solving the operator equation 
$$B_{\a,t} f  = g,$$
which is a well-posed problem.

\section{The Regularization} 

In view of the expression  (\ref{eq-5-0}),   for each $\a\in (0, 1)$ and $t\in (0, \t)$, we define the map 
$R_{t,\a}: L^2[0, \pi]\to L^2[0, \pi]$ as 
\beq\label{reg}   R_{t,\a} \psi = \sum_{n=1}^\infty \frac{ E_\a(-\l_n^2t^\a) } {E_\a(-\l_n^2\t^\a)}  \< \psi, \f_n\> \f_n,\q \psi \in L^2[0, \pi].\eeq
In view of (\ref{eq-6-1}), we see that $R_{t,\a}: L^2[0, \pi]\to L^2[0, \pi]$ is a well-defined bounded linear operator with 
$$\|R_{t,\a}\| 
 \leq \frac{C_2}{C_1}\Big(\frac{1+\t\a}{t^\a} \Big).$$
 Note that if $g = u_\a(\cdot, \t)$, then by (\ref{eq-5-0}), 
 $$R_{t,\a} g = u_\a(\cdot, t),\q 0<t<\t.$$
The following theorem shows that, for each $\a\in (0, 1)$,  the family $\{R_{t,\a}: 0<t<\t\}$  of operators defined above is a regularization family for the ill-posed inverse problem $P_0$.

\bt\label{Th-3} If $g=u_\a(\cdot, \t)$, then
$$\|u_\a(\cdot, 0)-u_\a(\cdot, t)\| \to 0\q\h{as}\q t\to 0.$$
\et 

\bpf
By the representation of $u_\a(\cdot, t)$ in  (\ref{eq-2}), we have 
$$u_\a(\cdot, t) - u_\a(\cdot, 0) = \sum_{n=1}^\infty [E_\a(-\l_n^2t^\a)-1] \<f_0, \f_n\>\f_n.$$
so that 
\beq\label{eq-10} \|u_\a(\cdot, t) - u_\a(\cdot, 0)\|^2  = \sum_{n=1}^\infty |E_\a(-\l_n^2t^\a)-1|^2| \<f_0, \f_n\>|^2 .\eeq
From the definition of $E_\a(\cdot)$, we have  
\beqarray
E(z)-1 & =&  \sum_{k=1}^\infty \frac{z^k}{\Gamma(\a k+1)}   =  z \sum_{k=1}^\infty \frac{z^{k-1}}{\Gamma(\a k+1)}  \\ 
&=&  z \sum_{k=0}^\infty \frac{z^k}{\Gamma(\a k+\a +1)}  = z E_{\a, \a+1}(z),
\eeqarray
where for $\a, \b>0$, $E_{\a, \b}(\cdot)$ is the {\it generalized Mittag-Leffler function} defined by 
$$E_{\a, \b}(z) = \sum_{k=0}^\infty \frac{z^k}{\Gamma(\a k+\b)} . $$
Hence, we have 
\beq\label{eq-12} E_{\a}(-\l_n^2t^\a)-1 = (-\l_n^2t^\a)  E_{\a, \a + 1}(-\l_n^2t^\a).\eeq
It is known that (cf. \cite{ARS}) that if  there exists $C>0$ such that for $0<\a<2$ and $\b\in \R$, 
$$|E_{\a, \b}(z)| \leq \frac{C}{1+|z|}$$
for all $z\in \C$ and for $\mu<|\arg(z)|\leq \pi$, where 
$$\frac{\p \a}{2}<\mu < \min\{\pi, \pi\a\}.$$
Now, taking  $z=-\l_n^2t^\a$ and $0<\a<1$, we have $\arg(z) = \pi$ and $\min\{\pi, \pi\a\} = \pi\a$.  Hence, in this case the required conditions on $\arg(z)$ is 
automatically satisfied. Thus, we have 
$$|E_{\a, \a+1}(-\l_n^2t^\a)|\leq \frac{C}{1+\l_n^2t^\a}.$$
Hence, 
\beq\label{eq-12-1}|E_{\a}(-\l_n^2t^\a)-1| = |-\l_n^2t^\a E_{\a, \a+1}(-\l_n^2t^\a)|\leq \frac{C\l_n^2t^\a}{1+\l_n^2t^\a}\leq C. \eeq
In particular, $|E_{\a}(-\l_n^2t^\a)-1|\to 0$ as $t\to 0$ for each $n\in \N$. 
Thus, 
$$|E_\a(-\l_n^2t^\a)-1|^2| \<f_0, \f_n\>|^2\to 0\h{ as }\, t\to 0\h{ for each }\, n\in \N$$
and
$$|E_\a(-\l_n^2t^\a)-1|^2|\<f_0, \f_n\>|^2  \leq  C^2 |\<f_0, \f_n\>|^2$$
with $\sum_{n=1}^\infty |\<f_0, \f_n\>|^2\leq \|f_0\|^2$. Hence, by the dominated convergence theorem, the relation (\ref{eq-10}) implies that 
$$\|u_\a(\cdot, t) - u_\a(\cdot, 0)\|\to 0\h{ as }\q t\to 0.$$
This completes the proof.
\epf 


\section{Error Estimate under Noisy Data}

If the data is noisy, say we have $\tilde g\in L^2[0, \pi]$ in place of $g$ such that 
$$\|g-\tilde g\| \leq \d$$
for some known noise level $\d>0$, then using the expression in (\ref{eq-5-0}), the corresponding solution at $t$ can be taken as 
$$\tilde u_\a(\cdot,t) = R_{\a,t}\tilde g,$$
where  $R_{\a,t}$ is defined as in (\ref{reg}). Thus,
\beq\label{noisy-sol} \tilde u_\a(\cdot,t)   = \sum_{n=1}^\infty \frac{ E_\a(-\l_n^2t^\a) } {E_\a(-\l_n^2\t^\a)}  \< \tilde g, \f_n\> \f_n.\eeq
Hence, we obtain 
$$\|u_\a(\cdot, t) - \tilde u_\a(\cdot,t)\|^2 = \sum_{n=1}^\infty \Big| \frac{ E_\a(-\l_n^2t^\a) } {E_\a(-\l_n^2\t^\a)}\Big|^2  |\< g- \tilde g, \f_n\>|^2 $$
so that using (\ref{eq-6-1}), 
$$\|u_\a(\cdot, t) - \tilde u_\a(\cdot,t)\|   \leq \frac{C_2}{C_1}\Big(\frac{1+\t^\a}{t^\a} \Big) \d .$$
Thus, we have proved:

\bt\label{Th-4} For $0<t<\t$, 
$$\|u_\a(\cdot, t) - \tilde u_\a(\cdot,t)\|   \leq \frac{C_2}{C_1}\Big(\frac{1+\t^\a}{t^\a} \Big) \d .$$
\et 

Thus we have proved the   following theorem.

\bt\label{Th-4-1} For $0<t<\t$, 
$$\|u_\a(\cdot, 0) - \tilde u_\a(\cdot,t)\|   \leq  \|u_\a(\cdot, 0) - u_\a(\cdot,t)\|  + \frac{C_2(1+\t^\a)}{C_1} \frac{ \d}{t^\a},$$
where (by Theorem \ref{Th-3}) $\|u_\a(\cdot, 0) - u_\a(\cdot,t)\|  \to 0$ as $t\to 0$.
\et 

\section{Error Estimates Under Source Conditions}

\noi
{\bf Assumption (A):} There exists an index function $\f: (0, \infty)\to [0, \infty)$ such that 
\beq\label{eq-13}\|u_\a(\cdot, t) - f_0\|  \leq c_0 \f(t)\eeq
for some $c_0>0$. 

\bt\label{Th-5}
Under the Assumption (A), let $\psi(t):= t\f(t)$ for $t>0$, and for $\d>0$, let  $t_\d:= \psi^{-1}(\d)$. Then 
$$\|f_0 - \tilde u_\a(\cdot,t_\d)\|  = O(\f(\psi^{-1}(\d)).$$
\et 

\bpf
By (\ref{eq-13})  and Theorem \ref{Th-4}, 
\beqarray \|f_0  - \tilde u_\a(\cdot,t)\|   
&\leq&  \|u_\a(\cdot, t) - u_\a(\cdot, 0)\|  +   \|u_\a(\cdot, t) - \tilde u_\a(\cdot,t)\|  \\ 
&\leq&  c_0 \f(t)   +  \frac{}{}  \frac{C_2(1+\t^\a)}{C_1} \frac{\d}{t^\a} 
\eeqarray
Note that  
$$\f(t) = \frac{\d}{t^\a} \iff \psi(t):= t^\a\f(t) = \d.$$
Hence, by choosing $t = t_\d:= \psi^{-1}(\d)$, we obtain 
$$\|f_0 - \tilde u_\a(\cdot,t)\|  =  O(\f(t_\d)) = O(\f(\psi^{-1}(\d)).$$
\epf

Now, we specify an index function $\f$ and a  source set  $M_\f$  such that (\ref{eq-13}) is satisfied whenever $f_0\in M_\f$.

Let $K_\a : L^2[0, \pi]\to L^2[0, \p]$ be defined by 
$$K_\a \f  = \sum_{n=1}^\infty E_\a(-\l_n^2\t^\a)  \<\f, \f_n\>\f_n,\q \f\in L^2[0, \pi].$$
Since $E_\a(-\l_n^2\t^\a) \in \R$ and  $E_\a(-\l_n^2\t^\a) \to 0$ as $n\to \infty$, it follows that (cf. \cite{Nair-fa})  $K_\a$ is a compact, self-adjoint operator and 
\beq\label{eq-source-eig} K_\a \f_n =  E(-\l_n^2\t^\a)\f_n\q\forall\, n\in \N.\eeq
In view of   (\ref{eq-3}), we see that $f_0:=u_\a(\cdot,0)$ is the solution of the compact operator equation 
$$K_\a f = g.$$  
Let 
\beq\label{eq-14} {\mathcal M}_{\a,\r}:=\{ K_\a u: \|u\| \leq \r\}.\eeq
and assume that 
$$f_0\in {\mathcal M}_{\a,\r}.$$
Then $f_0 = K_\a u$ for some  $u\in L^2[0, \pi]$ with $\|u\|\leq \r$, and  hence by (\ref{eq-source-eig}),  we have 
$$\<f_0, \f_n\> = \<K_\a u, \f_n\>  = \<u, K_\a \f_n\>  = E_\a(-\l_n^2\t^\a)\<u, \f_n\>\q\forall\, n\in \N.$$
Hence, from (\ref{eq-10}), we have 
\beqarray
 \|u_\a(\cdot, t) -  f_0 \|^2 & =&  \sum_{n=1}^\infty |E_\a(-\l_n^2t^\a)-1|^2| \<f_0, \f_n\>|^2\\ 
& =&  \sum_{n=1}^\infty |E_\a(-\l_n^2t^\a)-1|^2| |E_\a(-\l_n^2\t^\a)|^2  \<u, \f_n\>|^2
\eeqarray 
Now, (\ref{eq-12-1}) and Lemma \ref{lem-1} imply that 
$$|E_{\a}(-\l_n^2t^\a)-1|\, |E_\a(-\l_n^2\t^\a)|  \leq C_\a \frac{\l_n^2t^\a}{(1+\l_n^2t^\a) (1+\l_n^2\t^\a)},  $$
where  
\beq\label{eq-15} C_\a:=  \frac{CC_2}{\Gamma(1-\a)},\eeq
with $C$ and $C_2$ as in (\ref{eq-12}) and Lemma \ref{lem-1}, respectively.
Note that 
$$ \frac{\l_n^2t^\a}{(1+\l_n^2t^\a) (1+\l_n^2\t^\a)}  =  \frac{(\l_n^2t^\a)/(\l_n^2\t^\a)}{(1+\l_n^2t^\a)}  \frac{\l_n^2\t^\a} {(1+\l_n^2\t^\a)} \leq \frac{t^\a}{\t^\a} .  $$
Thus, we arrive at the estimate 
$$ \|u_\a(\cdot, t) -  f_0 \|^2   \leq C_\a^2 \Big(\frac{t^\a}{\t^\a}\Big)^2 \sum_{n=1}^\infty  | \<u, \f_n\>|^2.$$
Thus, we have proved the following theorem.

\bt\label{Th-est}
If $f_0\in {\mathcal M}_{\a,\r}$, then 
$$ \|u_\a(\cdot, t) -  f_0 \|   \leq \r C_\a   \frac{t^\a}{\t^\a},$$
where ${\mathcal M}_{\a,\r}$ and $C_\a$ are  as in (\ref{eq-14}) and (\ref{eq-15}), respectively. 
\et 

\brem\rm 
 Theorem \ref{Th-est} shows that  the function $\f$ defined by 
 $$\f(t) = t^\a,\q t>0,$$
 satisfies the Assumption (A) with $c_0 =  \r C_\a/\t^\a$. 
\erem

In view of Theorem \ref{Th-est} and Theorem \ref{Th-4-1}, if $f_0\in {\mathcal M}_{\a,\r}$, then we have 
\beqarray \|f_0  - \tilde u_\a(\cdot,t)\|   
&\leq&  \|u_\a(\cdot, t) - u_\a(\cdot, 0)\|  +   \|u_\a(\cdot, t) - \tilde u_\a(\cdot,t)\|  \\ 
&\leq&  \r C_\a   \frac{t^\a}{\t^\a}   +   \frac{C_2(1+\t^\a)}{C_1} \frac{\d}{t^\a}. 
\eeqarray
Now, 
$$\frac{t^\a}{\t^\a}  = \frac{\d}{t^\a}\iff t^{2\a} = \t^\a \d \iff t^\a = \sqrt{\t^\a\d} \iff \frac{t^\a}{\t^\a}  = \sqrt{\frac{\d}{\t^\a}}.$$
Thus, we have proved the following theorem.

\bt\label{Th-est-opt} 
If $f_0\in {\mathcal M}_{\a,\r}$ and $t_\d:= \sqrt\t  \d^{1/\a}$, then 
$$ \|\tilde u_\a(\cdot, t_\d) -  f_0 \|   \leq  \Big( \r C_\a     +  \frac{C_2(1+\t)}{C_1} \Big) \sqrt{\frac{\d}{\t^\a}} .  $$
where ${\mathcal M}_{\a,\r}$ and $C_\a$ are  as in (\ref{eq-14}) and (\ref{eq-15}), respectively. In particular, 
$$ \|\tilde u_\a(\cdot, t_\d) -  f_0 \| = O(\sqrt\d).$$
\et 

\brem\rm 
By the definition  of  ${\mathcal M}_{\a,\r}$ and from the standard regularization  theory  (cf. \cite{EHN, Nair-linop}), it follows that the estimate obtained in Theorem \ref{Th-est-opt} is optimal for the source set   ${\mathcal M}_{\a,\r}$.
\erem

\section{Numerical Illustrations }

In this section, we shall consider some numerical examples to illustrate the level of approximation of the regularized solutions.  

For numerical computations, we divide the given space domain $[0, \pi]$ into a finite number of equal subintervals with step size $h$ where $h = x_{i+1} - x_i, \; i = 0, 1, 2, \ldots, N-1$ with $x_0 = 0$, and $x_N = \pi$. All the simulations are carried out using MATLAB R2022a with step size $h = \pi/100$. 

We have observed that the problem of finding $u_{\alpha}(\cdot, t), \, 0 < t \leq \tau$ from the knowledge of $g := u_{\alpha}(\cdot, \tau)$ is a well-posed problem, and $\|u_{\alpha}(\cdot, t) - f_0\| \rightarrow 0$ as $t \rightarrow 0$ (see Theorem 3.1). 

For computational purpose, we take $\tau = 1$ and consider $g = u_{\alpha}(\cdot, \tau)$ obtained from 
(\ref{eq-2}) by taking $f_0(x) =  x (\pi - x ) e^{-x}, \; 0 \leq x \leq \pi$, and then compute $u_{\alpha}(\cdot,t)$ according to the formula (\ref{eq-5-0}). 

To approximate the integrals involved in the computation of $g$ and $u_{\alpha}(\cdot, t)$, we make use of the composite trapezoidal rule. For the illustration of the convergence $\|u_{\alpha}(\cdot, t) - f_0\|_{L^2} \rightarrow 0$ as $t \rightarrow 0$, we 
take $t_i = 10^{-(i+2)}$ for $ i = 1, 2, \ldots, 7$ and compute the exact expression $f_0(x)$ as well as $u_{\alpha}(x, t_i)$ and show in Figures $1$ to $4$ for various $\alpha \in \left\{0.2, \; 0.4, \;  0.6, \; 0.8 \right \}$.  

\begin{figure}[H]
  \centering
  \includegraphics[width=1.1\linewidth,height=0.43\textheight]{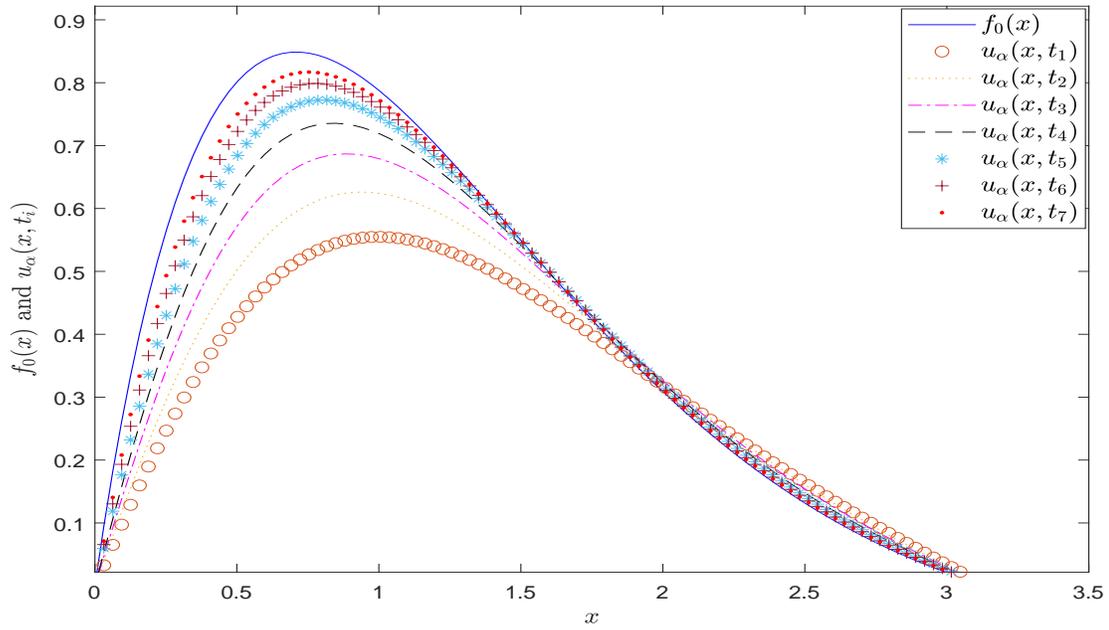}
  \caption{Solution profiles of $f_0(x)$ and $u_{\alpha}(x, t_{i})$ for $\alpha = 0.2$}
  \label{fig:EX5}
\end{figure}
\begin{figure}[H]
  \centering
  \includegraphics[width=1.1\linewidth,height=0.43\textheight]{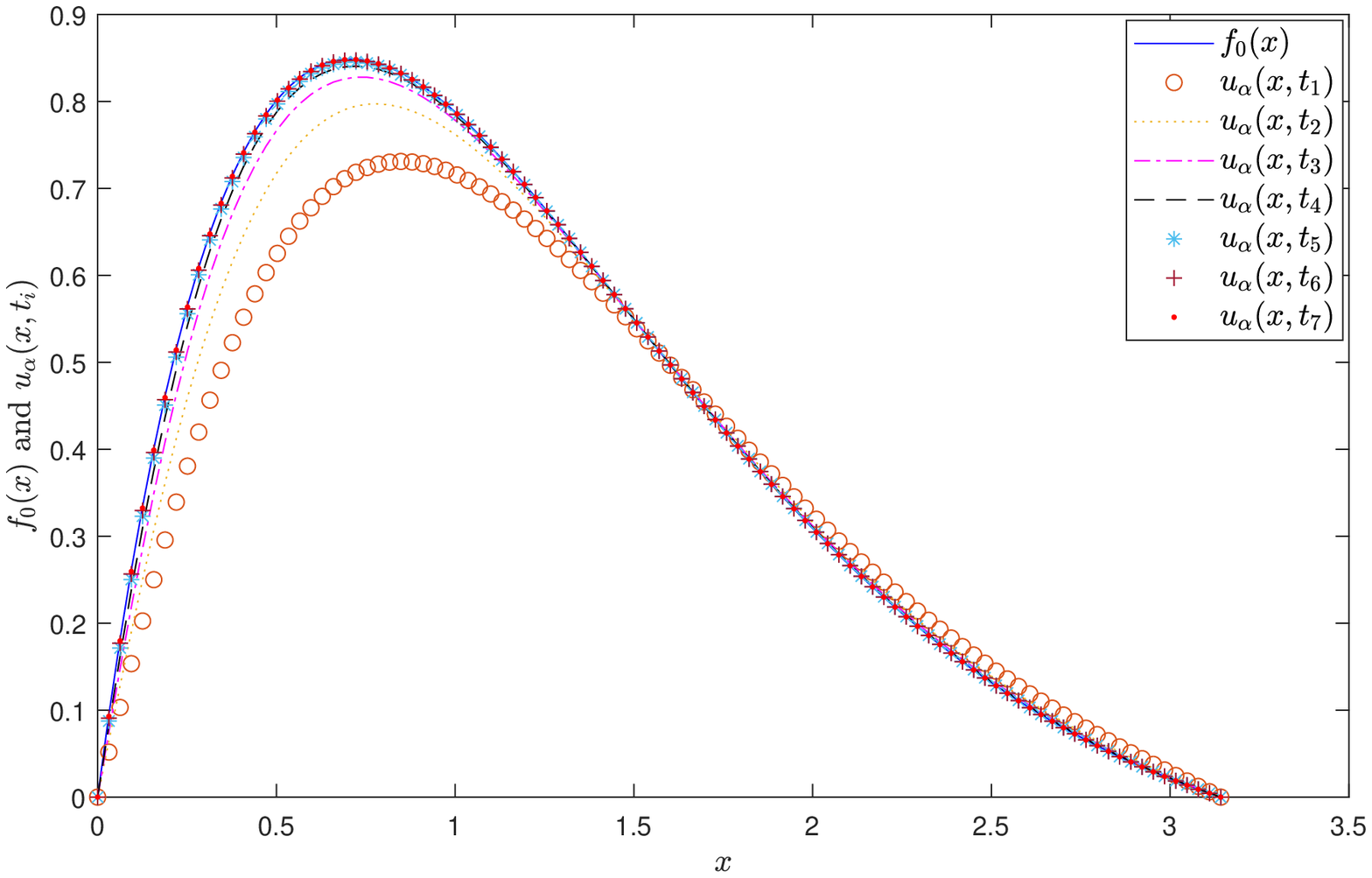}
  \caption{Solution profiles of $f_0(x)$ and $u_{\alpha}(x, t_{i})$ for $\alpha = 0.4$}
  \label{fig:EX6}
\end{figure}
\begin{figure}[H]
  \centering
  \includegraphics[width=1.1\linewidth,height=0.43\textheight]{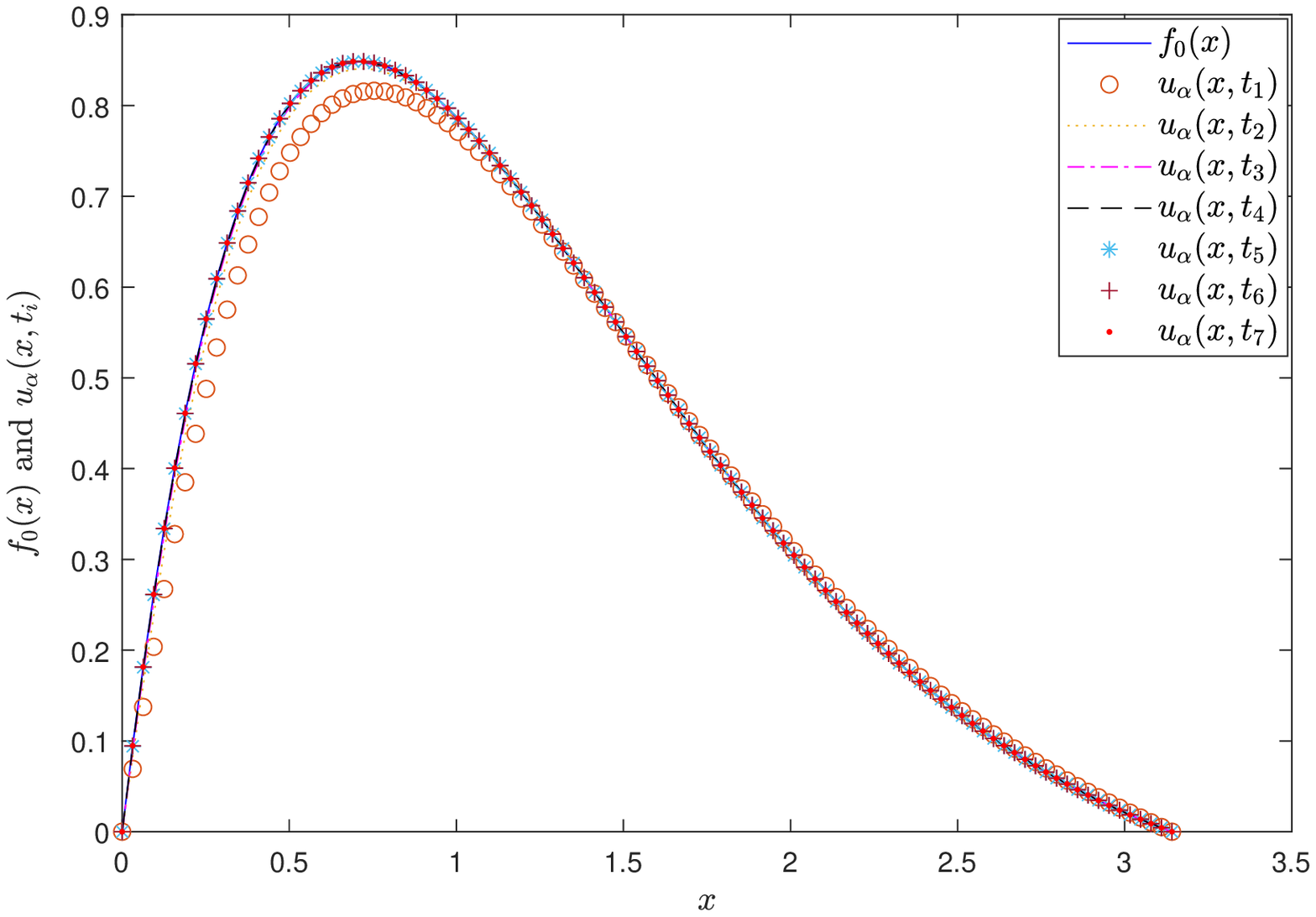}
  \caption{Solution profiles of $f_0(x)$ and $u_{\alpha}(x, t_{i})$ for $\alpha = 0.6$}
  \label{fig:EX7}
\end{figure}
\begin{figure}[H]
  \centering
  \includegraphics[width=1.1\linewidth,height=0.43\textheight]{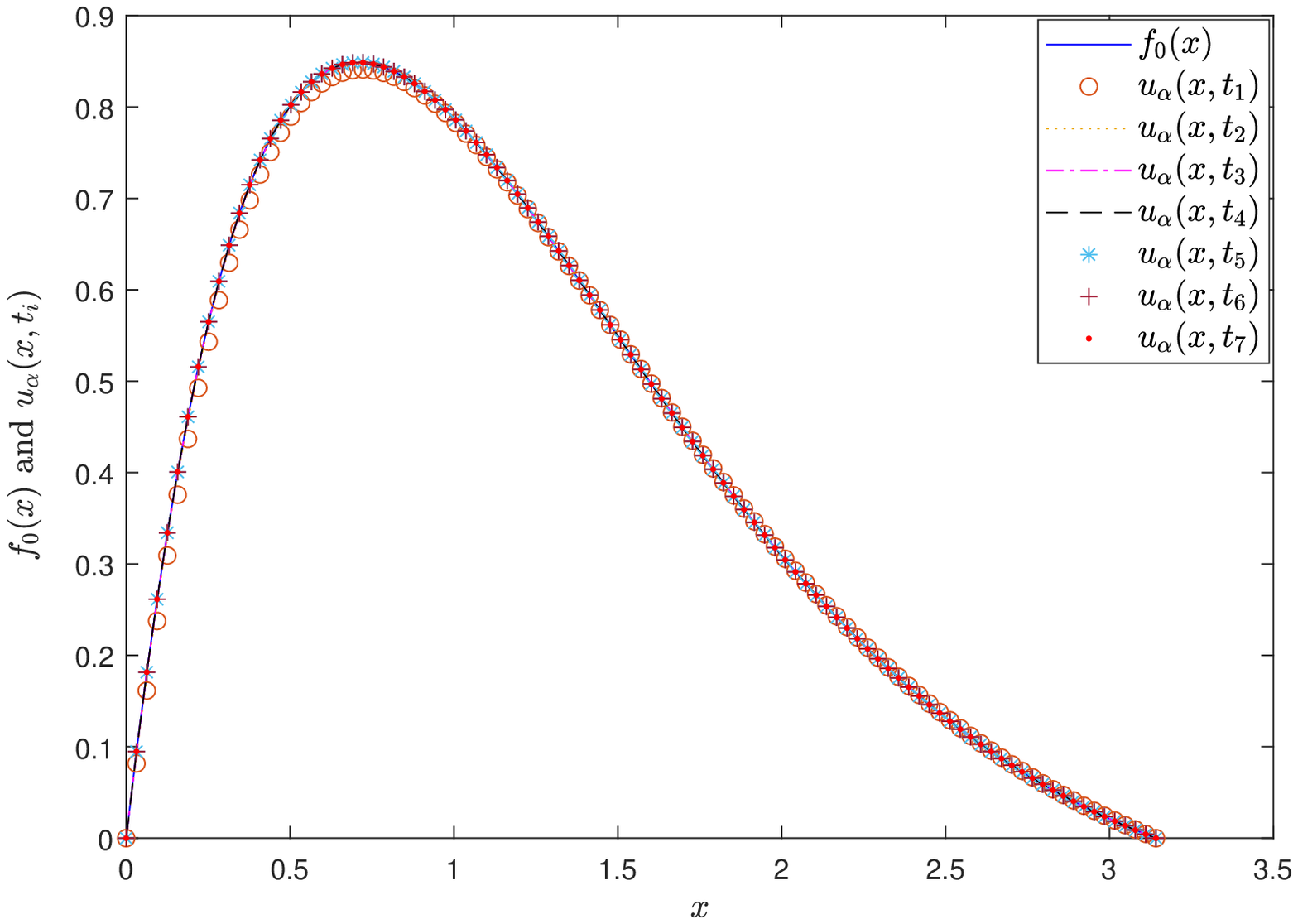}
  \caption{Solution profiles of $f_0(x)$ and $u_{\alpha}(x, t_{i})$ for $\alpha = 0.8$}
  \label{fig:EX8}
\end{figure}

\vsq 
In Table $1$, we show the error $\|u_{\alpha}(\cdot, t) - f_0 \|_{L^2}$ for  $t=t_i, \,  i = 1, 2, \ldots, 7$
for different values of $\alpha$. We observe that  when $t$  decreases, the error $\|u_{\alpha}(\cdot, t) - f_0 \|_{L^2}$ decreases. This validates our theoretical result in Theorem 4.1.

\begin{table}[H]
\begin{center}
\caption{$L^2$-errors between $f_0$ and $u_{\alpha}(\cdot, t_i)$ for $\alpha = 0.2, \; 0.4, \; 0.6$ and $0.8$}
\begin{tabular}{|c|c|c|c|c|} \hline
& $\alpha  = 0.2$ & $\alpha = 0.4$  & $\alpha  = 0.6$ & $\alpha  = 0.8$  \\ \cline{2-5}
$t_i$ \mbox{value} &   $\| u_{\alpha}(\cdot, t)-f_0\|_{L^2}$ & $\|f_0 - u_{\alpha}(\cdot, t)\|_{L^2}$ &  $\|u_{\alpha}(\cdot, t)-f_0\|_{L^2}$  &  $\| u_{\alpha}(\cdot, t) - f_0\|_{L^2}$  \\ \hline
$10^{-3}$ &  3.2829(-1) & 1.4958(-1) & 5.2366(-2)& 1.5080(-2) \\ 
$10^{-4}$ &  2.5593(-1) & 7.5081(-2) & 1.5476(-2)& 2.6205(-3) \\ 
$10^{-5}$ &  1.9370(-1) & 3.4924(-2) & 4.2280(-3)& 4.2730(-4) \\ 
$10^{-6}$ &  1.4253(-1) & 1.5401(-2) & 1.1058(-3)& 6.8572(-5) \\ 
$10^{-7}$ &  1.0225(-1) & 6.5509(-3) & 2.8063(-4)& 2.6230(-5) \\ 
$10^{-8}$ &  7.1751(-2) & 2.7197(-3) & 7.1304(-5)& 2.5992(-5) \\ 
$10^{-9}$ &  4.9403(-2) & 1.1102(-3) & 2.8563(-5)& 2.6255(-5) \\ 
\hline
\end{tabular}
\end{center}
\end{table}
\noindent

For the illustration of the case with  the noisy data, we take 
$
{\tilde g}(x) = g(x) + \frac{\delta}{2} 
$
with some noise level $\delta > 0$. Note that $\|{\tilde g} - g\|_{L^2} \leq \delta$.
Now we take 
$t_{\delta} := \sqrt{\tau} \, \delta^{1/\alpha}$ as per Theorem 6.3 and  compute $\tilde{u}_{\alpha}(x, t_{\delta})$ using  the formula (\ref{noisy-sol}) for  several values of $\delta$ and for some values of  $\alpha$.  

Taking $\delta_i = 10^{-(i+2)}, i = 1, 2, \ldots, 7$  we  compute the exact expression $f_0$ as well as $\tilde{u}_{\alpha}(x, t_{\delta_i})$ and show them in Figures $5$ to $8$ for various $\alpha \in \left\{ 0.2, \; 0.4, \; 0.6, \; 0.8 \right \}$. 

\begin{figure}[H]
  \centering
  \includegraphics[width=1.1\linewidth,height=0.42\textheight]{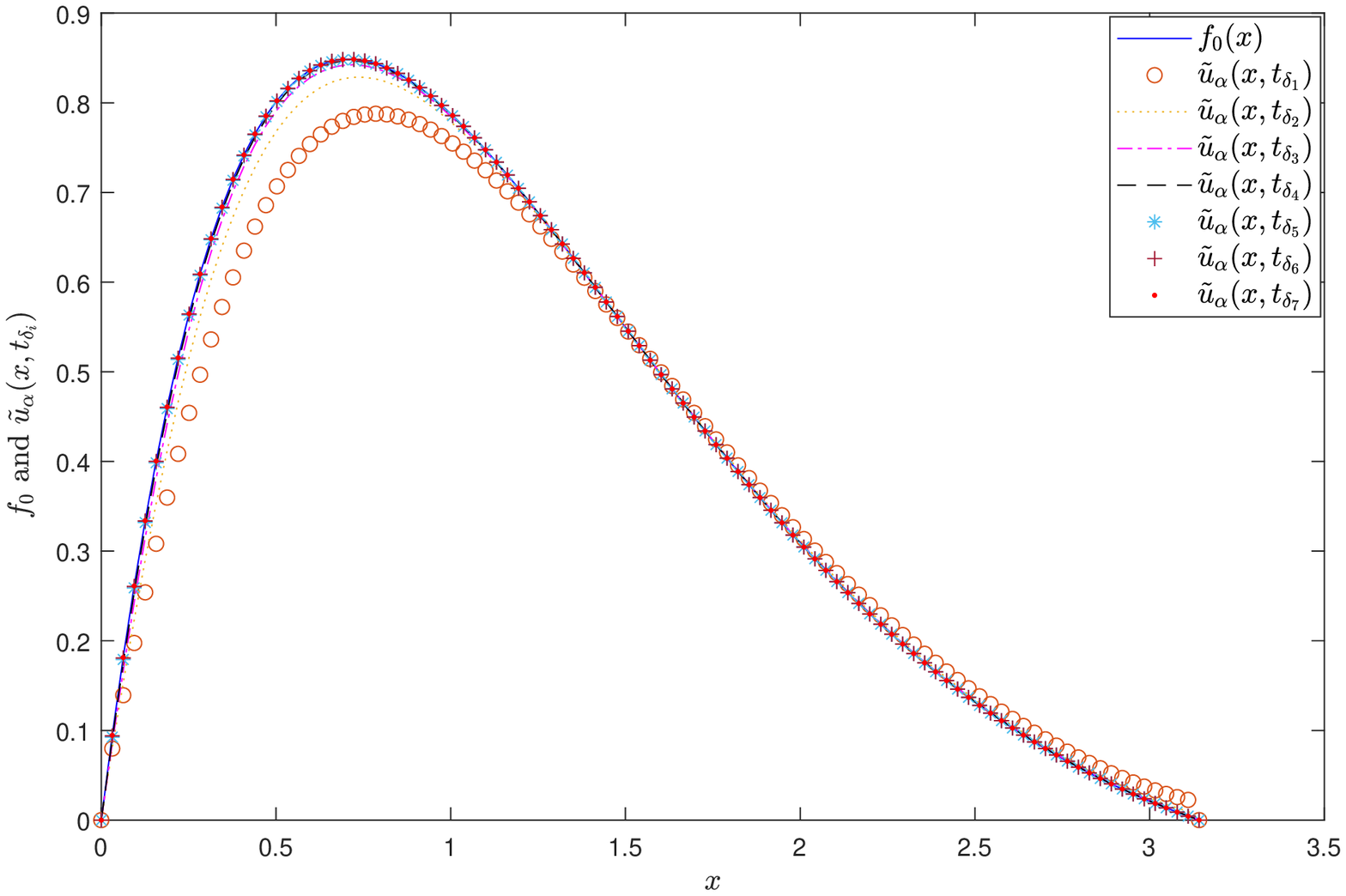}
  \caption{Solution profiles of $f_0(x)$ and $\tilde{u}_{\alpha}(x, t_{\delta_i})$ for $\alpha = 0.2$}
  \label{fig:EX9}
\end{figure}
\begin{figure}[H]
  \centering
  \includegraphics[width=1.1\linewidth,height=0.43\textheight]{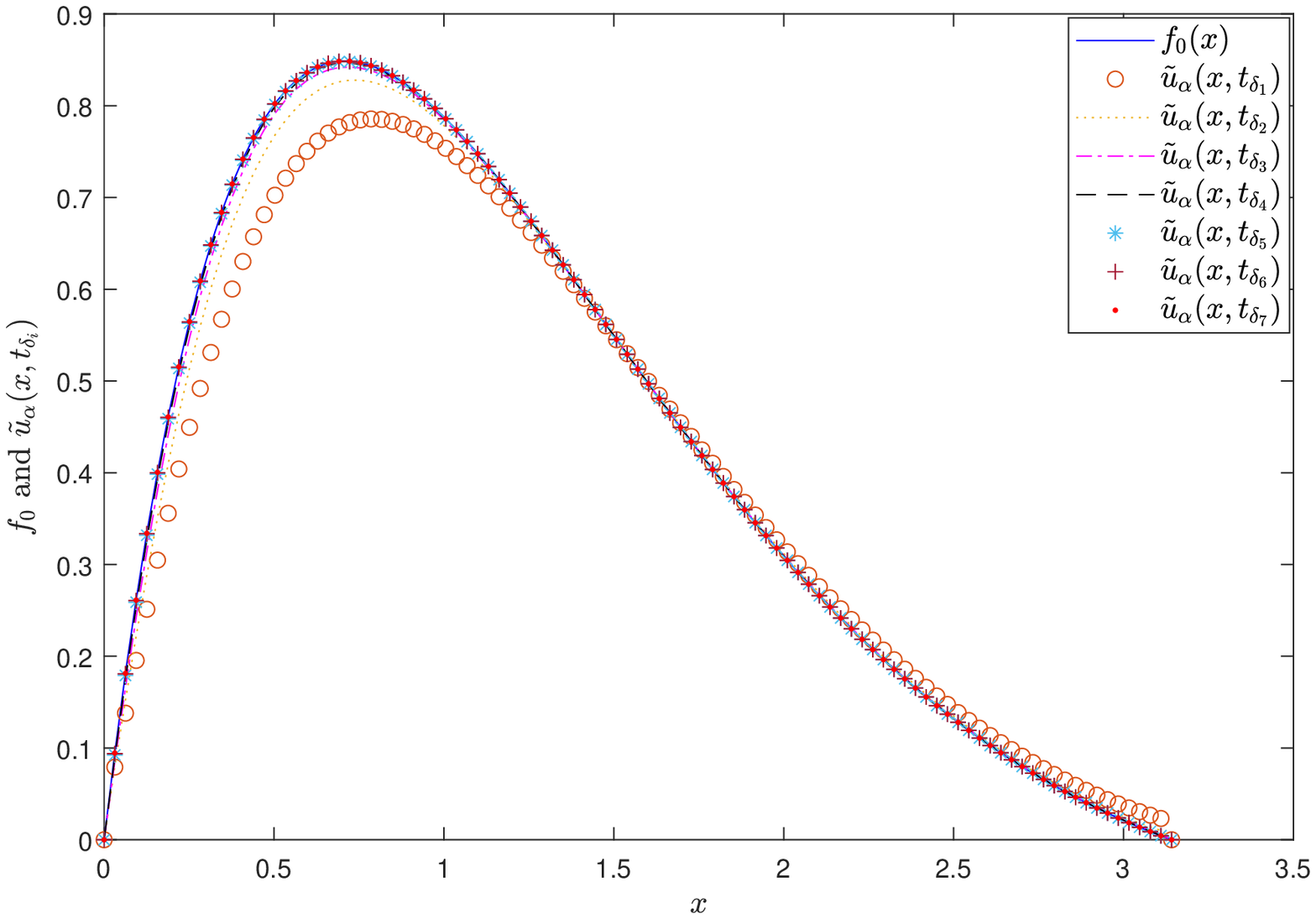}
  \caption{Solution profiles of $f_0(x)$ and $\tilde{u}_{\alpha}(x, t_{\delta_i})$ for $\alpha = 0.4$}
  \label{fig:EX10}
\end{figure}
\begin{figure}[H]
  \centering
  \includegraphics[width=1.1\linewidth,height=0.43\textheight]{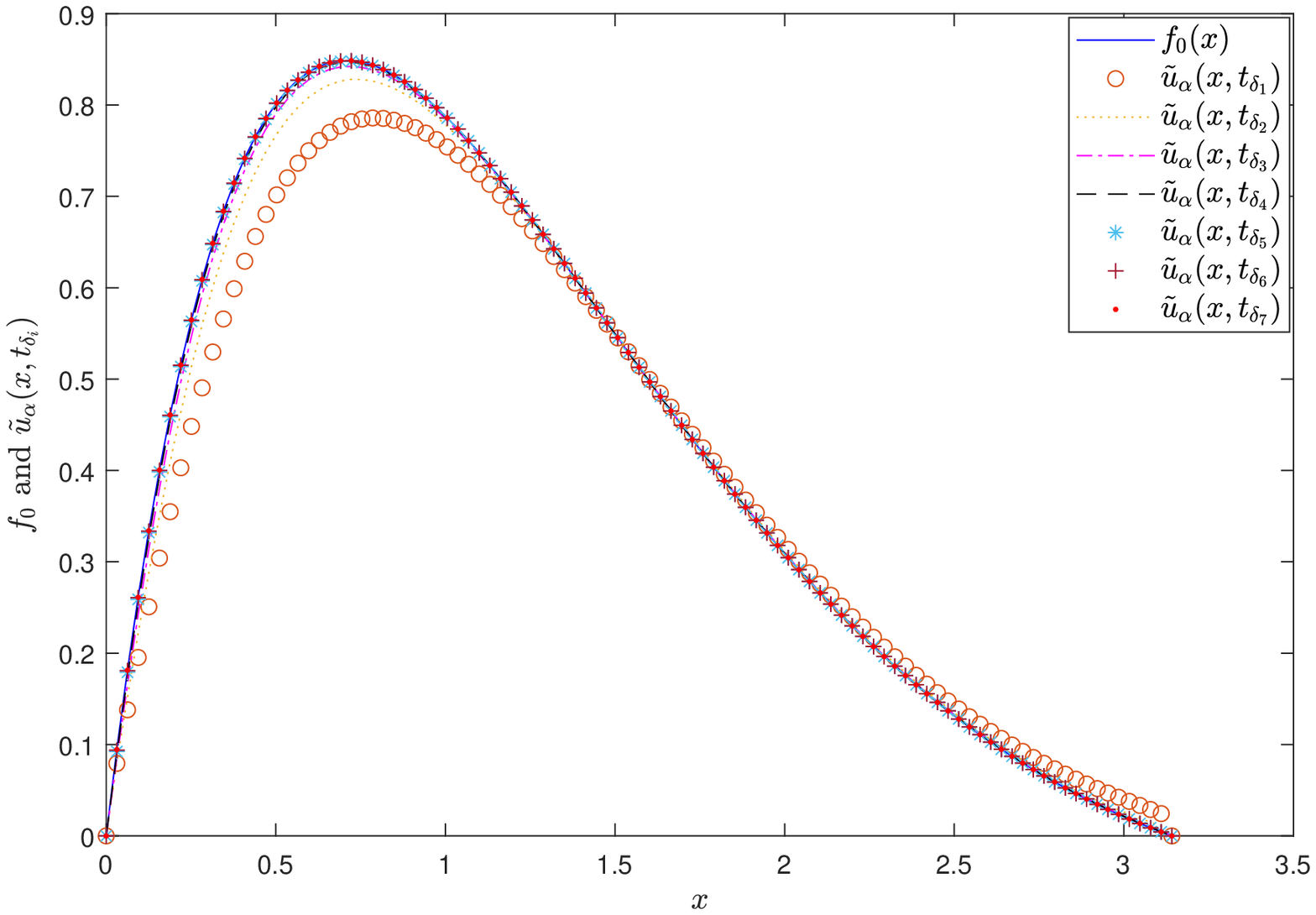}
  \caption{Solution profiles of $f_0(x)$ and $\tilde{u}_{\alpha}(x, t_{\delta_i})$ for $\alpha = 0.6$}
  \label{fig:EX11}
\end{figure}
\begin{figure}[H]
  \centering
  \includegraphics[width=1.1\linewidth,height=0.43\textheight]{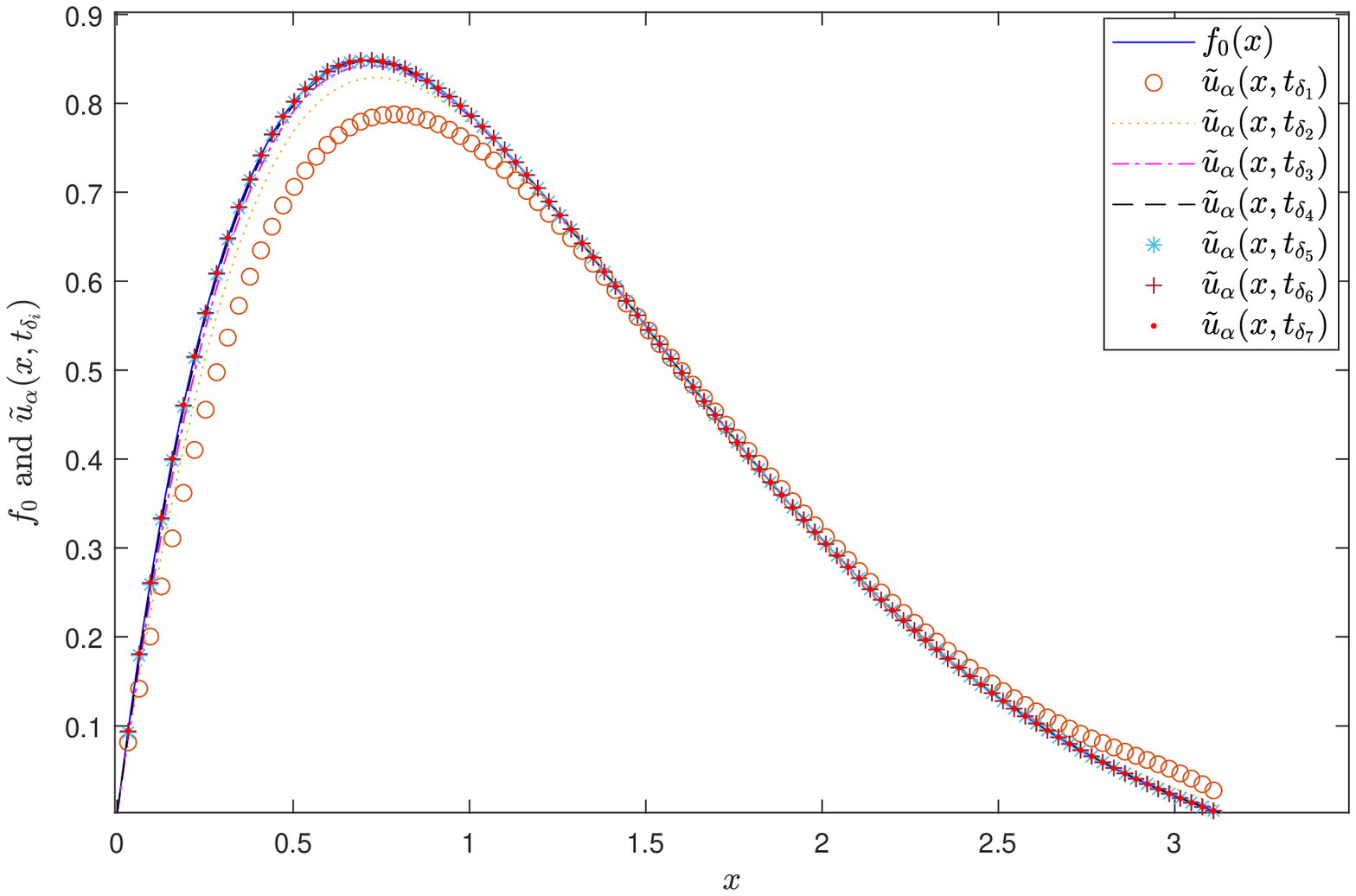}
  \caption{Solution profiles of $f_0(x)$ and $\tilde{u}_{\alpha}(x, t_{\delta_i})$ for $\alpha = 0.8$}
  \label{fig:EX12}
\end{figure}

\vsq 

In Table $2$, we show the error $\|\tilde{u}_{\alpha}(\cdot, t_{\delta}) - f_0 \|_{L^2}$ 
for different values of $\alpha$ and $\delta_i, i = 1, 2, \ldots, 7$. Note that when $\delta$ value decreases, the error $\|\tilde{u}_{\alpha}(\cdot, t_{\delta}) - f_0 \|_{L^2}$ decreases. This validates our theoretical result in Theorem 6.3. 

\begin{table}[H]
\begin{center}
\caption{$L^2$-errors between $f_0$ and $\tilde{u}_{\alpha}(\cdot, t_{\delta_i})$ for $\alpha = 0.2, \; 0.4, \; 0.6$ and $0.8$}
\begin{tabular}{|c|c|c|c|c|} \hline
& $\alpha  = 0.2$ & $\alpha = 0.4$  & $\alpha  = 0.6$ & $\alpha = 0.8$  \\ \cline{2-5}
$\delta_i$ \mbox{value} &   $\| \tilde{u}_{\alpha}(\cdot, t_{\delta}) - f_0 \|_{L^2}$ & $\| \tilde{u}_{\alpha}(\cdot, t_{\delta}) - f_0 \|_{L^2}$ &  $\|f_0 - \tilde{u}_{\alpha}(\cdot, t_{\delta})\|_{L^2}$ &  $\| \tilde{u}_{\alpha}(\cdot, t_{\delta})-f_0 \|_{L^2}$  \\ \hline
 $10^{-3}$ & 8.2528(-2)  & 8.6054(-2) & 8.6816(-2) & 8.3284(-2)  \\ 
$10^{-4}$ &  3.2857(-2) & 3.4183(-2) & 3.4288(-2) & 3.2853(-2)\\ 
$10^{-5}$ &  1.1881(-2) & 1.2337(-2) & 1.2327(-2) & 1.1816(-2)\\ 
$10^{-6}$ &  4.0572(-3) & 4.2074(-3) & 4.1935(-3) & 4.0215(-3)\\ 
$10^{-7}$ &  1.3385(-3) & 1.3870(-3) & 1.3802(-3) & 1.3240(-3) \\ 
$10^{-8}$ &  4.3110(-4) & 4.4651(-4) & 4.4386(-4) & 4.2573(-4)\\ 
$10^{-9}$ &  1.3645(-4) & 1.4123(-4) & 1.4027(-4) & 1.3453(-4) \\ 
\hline
\end{tabular}
\end{center}
\end{table}


\vsq
\noindent 
{\bf Acknowledgement:} \, The first author M. Thamban Nair gratefully acknowledges the support received from BITS Pilani, K.K. Birla Goa Campus, where he is a Visiting Professor from August 1, 2023 after superannuation from I.I.T. Madras, Chennai.

\end{document}